\documentstyle[12pt]{article}

\newcommand{\be}{\begin{equation}}
\newcommand{\ee}{\end{equation}}
\newcommand{\zz} {\overline{z}}
\newcommand{\nn} {\noindent}
\begin{document}
\title{Bargmann representation for some deformed harmonic oscillators\\
with non-Fock representation \footnote{ presented at the Symposium in honour of Jiri Patera and Pavel Winternitz for their 60th birthday , Algebraic Methods and Theoretical Physics, January 9-11, 1997, Centre de recherches math\' ematiques, Universit\' e de Montr\' eal.}}
\author{Mich\`{e}le IRAC-ASTAUD and Guy RIDEAU\\ 
Laboratoire de Physique Th\'{e}orique de la mati\`ere condens\'ee\\
Universit\'{e} Paris VII\\2 place Jussieu F-75251 Paris Cedex 05, FRANCE}
\date{}
\maketitle
Preprint  1997. \\
to appear in the Proceedings of the Symposium.
\begin{abstract}
We prove that Bargmann representations exist for some deformed harmonic oscillators that admit non-Fock representations. In specific  cases, we explicitly obtain the resolution of the identity   in terms of a true integral on the complex plane. We prove on explicit examples that Bargmann representations cannot always be found, particularly when the coherent states do not exist in the whole complex plane.
\end{abstract}

\section{Introduction}
We define a deformed harmonic oscillator  by the following algebraic relations between three operators $a,a^\dagger$ and $N$:

\begin{equation}[a,N] = a ,\quad [a^\dagger ,N] =
-a^\dagger \label{a1}
\end{equation}
\nn and
\be
a^\dagger a = \psi(N) , \quad aa^\dagger =\psi(N+1)
\label{a3}
\ee
\nn where $\psi$ is some real function and $a^\dagger$ has to be the adjoint of $a$ and $N$ is self-adjoint.

\nn The algebra defined in (\ref{a1}), (\ref{a3}) was studied in \cite {nous1} and related algebras, under different forms, were considered in \cite{nous2}\cite {quesne}.

\nn In the  case of the usual harmonic oscillator, $\psi(N)=N$, a particularly interesting representation of (\ref{a1}), (\ref{a3}) is the Bargmann representation, realized in a Hilbert space of the entire functions of order $\leq 1$ with a reproducing kernel.

\nn The main purpose of this paper is to study if similar representations can be obtained when $\psi(N)$ is a strictly positive function. In this case, we are faced with non-Fock representations since in any representation the spectrum of $N$ has no lower and upper bounds. In the next section, we give a brief description of these representations and we discuss the existence of the corresponding coherent vectors. In section 3, we discuss in some detail the possibility of a Bargmann representation when the function $\psi$ goes from $0$ to $\infty$ when $x$ goes from $-\infty$ to $+\infty$. In section 4, we prove the existence of a Bargmann representation in the particular case of a q-oscillator. In section 5, we investigate other examples proving that a Bargmann representation does not always exist. In section 6, we briefly discuss few examples  where $\lim _{x \rightarrow -\infty}\psi(x)> 0$ and/or $\lim _{x \rightarrow +\infty}\psi(x)< \infty$ . In such cases, though  cohere!
 nt!
 sta
tes exist but associated with a bounded part of the  complex plane, we can prove the impossibility of getting a Bargmann representation.

\nn It is not irrelevant to precise that the choice, although usual, of the letter $N$ is rather misleading as apparently referring (in the Fock representation) to the number of particles. The eigenvalues of $N$  label the levels of the energy operator that is not a priori given.

\section{Representations}

Let us assume $N$ has an eigenvector $\mid 0>$, with eigenvalue $\mu$.
The inequivalent representations are labelled by  the decimal part of $\mu$. The Hilbert space ${\cal H}$ of the representation is spanned by the normalized eigenvectors  $\mid n >$  given by 
\be
\mid n> =\left\{
\begin{array}{ll}
\lambda _n a^{\dagger n}\mid 0>,& \quad n\in Z+\\
\lambda _n a^{-n}\mid 0>,& \quad n\in Z-
\end{array}
\right.
\ee
\nn with 
\be
\lambda_n^{-2} =\psi (\mu +n)!= \left\{\begin{array}{ll}
\prod_{i=1}^n \psi(\mu +i),& \quad n \in Z^+\\
 \prod_{i=0}^{n+1} \psi(\mu +i),& \quad n \in Z^-
\end{array}
\right.
\ee
 \nn We have
\be
 \left\{\begin{array} {ll}
a^\dagger \mid n >  =& (\psi (\mu +n+1)^{1/2} \mid n+1 > \\ 
a\mid n >  =& (\psi (\mu + n))^{1/2} \mid n-1 > ,\quad n \in Z \\ 
N \mid n >  =&  (\mu +n )\quad \mid n > 
 \end{array}
\right.
\label{repre}
\ee

\nn The first step towards a Bargmann representation requires to study the coherent vectors, i.e, the eigenvectors of the annihilation operator $a$. Let us denote by $\mid z>$ the vector
\be
 a \mid z> = z \mid z>
\ee
\nn $\mid z>$ is given up to a constant of normalization by
\be
\mid z> =\sum_{n=-1}^{-\infty}z^{n}(\psi(\mu +n)!)^{1/2}\mid n>+\sum_{n=0}^{\infty}z^n (\psi(\mu +n)!)^{-1/2}\mid n>
\label{z}
\ee
\nn with the convention $\psi(\mu )! =1$.

\nn The vector $\mid z>$ belongs to ${\cal H}$ only if the two series in the right hand side of (\ref{z}) are simultaneously  norm convergent. This implies that :
\be
\mid z \mid <\lim_{p \rightarrow +\infty} \psi(p)^{1/2}\equiv r_2
\ee
\nn and
\be
\mid z \mid >\lim_{p \rightarrow -\infty} \psi(p)^{1/2}\equiv r_1
\ee
\nn So the coherent states exist only if $\psi$ is such that $r_1< r_2$
and their domain of existence is  a ring defined by $r_1<\mid z\mid < r_2$ which can be extended to the complex plane without the origin when $r_1=0$ and $r_2=\infty$.

\nn When $r_1$ is larger than $r_2$, the annihilation operator $a$ has no eigenstates  but then the creation operator $a^\dagger$ has eigenstates. Both situations are analogous and in the following we restrict to the case where $r_1< r_2$. 

\nn Although $\mu$ is a significant quantity for labelling inequivalent representations, it does not play a  part in the present problem. 
So we simplify the notation in assuming $\mu =0$ from now on. Indeed, this is equivalent to substitute $N-\mu$ to $N$ and $\psi_{\mu}(N)=\psi (\mu + N)$ to $\psi (N)$.

\nn Moreover we assume from now on that, unless otherwise specified, the function $\psi(N)$ is such that $r_1=0$ and $r_2=\infty$; this means that the coherent vectors exist for all complex $z$, except $z=0$.

\nn {\bf Remark :} ${\cal H}$ can be seen as the space of the functions $f(\theta ),0\leq \theta \leq 2\pi$ on the unit circle such that 
\be
\int_0^{2\pi} \mid f(\theta )\mid ^2 < \infty
\ee
\nn To the basis vectors $\mid n>$ correspond the functions $\exp (in\theta ) , n \in Z $ and the operators $a, a^\dagger $ and $N$ read :
\be
\left\{\begin{array} {ll}
a =& \exp (-i \theta)\psi (-i \frac{d}{d\theta} +\mu )^{1/2}\\
a^\dagger =& \psi (-i \frac{d}{d\theta} +\mu )^{1/2} \exp (i \theta)\\
N=& -i \frac{d}{d\theta}
 \end{array}
\right.
\ee
\nn which make sense since $\psi(x)$ is strictly positive. From this point of view, we can speak of a deformed rotator, in the same line of thought as in \cite{kowalski}.

\section{Towards a Bargmann representation}

 \nn Let us stress that in this paper we look for a Bargmann representation in the strict sense where the integrals involved in the scalar product, or equivalently in the resolution of identity are Riemann integrals. 

\nn So we are looking  for a positive real function $F$ , that we will call the weight function, such as we  have the following resolution of the identity :
\be
\int F(z \zz ) \mid \zz ><\zz \mid dz d\zz =1
\label{1}
\ee
\nn where the integration is extended to the whole complex plane.

\nn Then to a vector $\mid f>$ in the representation space $\cal{H}$
 \be
\mid f>= \sum _{-\infty}^{+\infty}f_n \mid n> , \quad \sum _{-\infty}^{+\infty} \mid f_n \mid^2< \infty
\ee
\nn should correspond the function of $z$ :
\be
f(z) \equiv < \zz \mid f > = \sum _{n \geq 0} z^n f_n (\psi (n)!)^{-1/2} + \sum _{n < 0} z^n f_n (\psi (n)!)^{1/2}
\label{f}
\ee
\nn In particular, to the basis vectors $\mid n>$ correspond  : 
\be
< \zz \mid n > = \left\{ \begin{array}{ll} 
z^n  (\psi (n)!)^{-1/2} , \quad &n \geq 0\\
z^n (\psi (n)!)^{1/2}, \quad &n <0
\end{array}
\right.
\label{n}
\ee
\nn From formula (\ref{1}) the scalar product in $\cal{H}$ reads :
\be
<g\mid f> = \int F(z \zz )\overline{g(z)}f(z) \mid dz d\zz 
\ee
and any operator $A$ is represented by the kernel $A(\zeta,\zz) = <\overline {\zeta}\mid A\mid \zz>$ such as~:
\be
 Af(\zeta ) \equiv <\overline{\zeta }\mid Af> = \int F(z \zz )A(\zeta, z) f(z) dz d\zz 
\ee
\nn Formula (\ref{f}) defines 
 a set $\cal{S}$ of holomorphic functions on the complex plane without the origin, the analytical properties of which are strongly depending on the function $G(x)$ defined by :
\be
G(x) = \sum _{n \geq 0} x^n  (\psi (n)!)^{-1} + \sum _{n < 0} x^n  \psi (n)!
\label{G}
\ee
\nn Indeed from the Schwarz inequality, we get: 
\be
\mid f(z) \mid \leq \mid f\mid ^{1/2}  G(z \overline{z}) ^{1/2}
\label{<}
\ee
\nn so that $G(x)^{1/2}$ controls the growth of the functions in ${\cal S}$ at infinity and near the origin.

\nn Moreover,  $G(\zeta z)=<\zz \mid \zeta >$ corresponds in ${\cal S}$  to the coherent state $\mid \zeta >$ ; and finally, if (\ref{1}) should be true, $G(\zeta \zz) $ is the reproducing kernel since we would have :
 \be
  f(\zeta) =\int F(z \zz ) < \overline{\zeta}\mid \zz > f(z) dz d\zz 
\ee
 \nn Taking (\ref{repre}) into account, we deduce that in our Bargmann representation, $a^\dagger$ would be the multiplication by $z$,  $a$ the operator $z^{-1}\psi(zd/dz)$ and $N$ the operator $zd/dz$. As $G(\zeta z)$  corresponds  to the coherent state $\mid \zeta >$, we obtain for $G$ the following functional equation
\be
x G(x) = \psi(x \frac{d}{dx}) G(x)
\label{eqG}
\ee
which could be obtained directly from the expansion (\ref{G}).

\nn From this set of inferences, we get the following necessary conditions to be verified by $F(x)$
 :
\be
1)\quad M(n)=\int_0^{\infty}  F(x) x^n dx = \left\{ \begin{array}{ll}
\psi(n)!, & n\geq 0\\
(\psi(n)!)^{-1}, & n <0
\end{array}
\right.
\label{moments}
\ee
\nn and
\be
2) \quad \int F(z \zz )  z f(z) \overline{g(z)} dz d\zz = \int F(z \zz ) f(z) \overline{\frac{1}{z}\psi(z \frac{d}{dz}) g(z)} dz d\zz
\label{adjoint}
\ee
\nn The first condition derives from (\ref{n}) and the second expresses the mutual adjointness of $a$ and $a^\dagger$. But Formulas (\ref{moments}) imply that $F(x)$ decreases fast  at the origin  and at infinity.  Then, at least if $f(z),g(z)$ are finite linear combinations of powers of $z$, we obtain : 
\be
 \int F(z \zz ) f(z) \overline{\frac{1}{z}\psi(z \frac{d}{dz}) g(z)} dz d\zz =
\int \left(\psi( -\zz \frac{d}{d\zz} -1) \frac{F(z \zz )}{z\zz} zf(z)\right)\overline{ g(z)} dz d\zz
\ee
\nn Using (\ref{adjoint}), we obtain the following functional equation for $F(x)$ :
\be
\psi(-x \frac{d}{dx}-1 ) \frac {F(x)}{x} = F(x)
\ee
\nn or equivalently
\be
x F(x) = \psi(-x \frac{d}{dx}) F(x)
\label{FF}
\ee
\nn for which we are looking for solutions fast decreasing at infinity and at the origin.

\nn Let us remark that (\ref{moments}) implies the following recursive relation between the momentums $M(n)$ of $F(x)$ :
\be
M(n+1)=\psi (n) M(n)
\label{M}
\ee
\nn Furthermore, according to the behaviour of $F(x)$ at infinity and at the origin, we can define its Mellin transform $\hat {F}(\rho)$ :
\be
\hat {F}(\rho)=\int_0^{\infty} F(x) x^{\rho -1} dx
\label{mel}
\ee
\nn which is meaningful for any $\rho$. According to (\ref{FF}), we must have :
\be
\hat {F}(\rho+1)= \psi(\rho)\hat {F}(\rho)
\label{mel1}
\ee
\nn As $\hat{F}(n+1)/\hat{F}(1) = M(n)$, $\hat{F}(\rho)$ is an interpolation of $M(n)$ which is univoquely determined by the condition of the mutual adjointness of $a$ and $a^\dagger$. It is worthwhile to note the following simple formula 
\be
G(x) = \hat{F}(1)\sum _{-\infty}^{+\infty}\frac{x^n}{\hat{F}(n+1)}
\label{GF}
\ee
 \nn According to our hypothesis, Equation (\ref{FF}) is not a differential equation. Indeed in this case $\psi (x)$ would be a polynomial. This is impossible for if the higher term is an odd power of $x$, $\psi$ cannot be a strictly positive function and if the highest term is an even power of $x$, coherent vectors does not exist. This implies that the general solution of (\ref{FF}) may depend from an arbitrary function. The examples of the next sections will be illuminating in this respect.

\section{The "q-oscillator"}

\nn In this section, we assume $\psi_{\lambda , q}(N)= \lambda q^{-N}, q \leq 1 , \lambda >0$. When $\lambda = (q^{-1}-q)^{-1}$, $a$ and $a^\dagger$ verify the well-known q-oscillator commutation relation~:
\be
a a^\dagger - q a^\dagger a = q^{-N}
\ee
\nn so that the function $(q^{-1}-q)^{-1} q^{-N}$ characterizes a particular non-Fock representation of the q-oscillator. Let us remark that we also have $a a^\dagger = q^{-1}a^\dagger a$, which corresponds to the complex q-plane.

\nn We  easily get :
\be
G_{\lambda , q}(x) = \sum_{-\infty}^{+\infty}\left(\frac{x}{\lambda}\right)^n q^{n(n+1)/2}
\ee
\nn which  obviously verifies the functional equation 
\be
x G_{\lambda , q}(x) = \lambda q^{-x\frac{d}{dx}}G_{\lambda , q}(x) = \lambda G_{\lambda , q}(q^{-1}x)
\label{Gq}
\ee
\nn Let us remark that 
\be
G^0_{\lambda , q}(x) = \exp\left(-\frac{\ln^2(x\lambda ^{-1})}{2\ln(q)}-\frac{\ln(x\lambda ^{-1})}{2}\right)
\ee
\nn verifies (\ref{Gq}). The general solution takes the form 
\be
G_{\lambda , q}(x)= q^{-1/8} G^0_{\lambda , q}(x) \sum _{-\infty}^{+\infty} \exp \left( \frac{1}{2} \ln (q)\left( n+  \frac{1}{2} +\frac{\ln (x\lambda^{-1})}{\ln q }\right)^2 \right)
\ee
\nn where the series in the right hand side is a strictly positive bounded periodic function of $\ln (x\lambda^{-1})/\ln q$ of period 1. So $G_{\lambda , q}(x)$ is increasing at infinity faster than any positive power of $x$ and at the origin faster than any negative power of $x$. Moreover, its growth is always less than  exponential of any power of $x$. Therefore, due to (\ref{<}), this type of growth characterizes the function $f(z)$ of the set $\cal{S}$.

\nn We solve now the functional equation (\ref{FF}):
\be
x F_{\lambda , q}(x) = \lambda F_{\lambda , q}(qx)
\label{Fq}
\ee
\nn The Mellin transform $\hat{F}_{\lambda , q}(\rho)$ verifies :
\be
\hat{F}_{\lambda , q}(\rho+1)= \lambda q^{-\rho}\hat{F}_{\lambda , q}(\rho) 
\ee
\nn with the particular solution :
\be
\hat{F}^0_{\lambda , q}(\rho)= \exp\left( \rho \\ln\lambda -\frac{1}{2}(\rho^2-\rho)\ln q \right)
\ee
\nn The inverse Mellin transform reads :
\be
F^0_{\lambda , q}(x) = \exp\left(\frac{\ln^2(x\lambda ^{-1})}{2\ln (q)}-\frac{\ln(x\lambda ^{-1})}{2}\right)
\label{F0}
\ee
\nn Therefore the general solution of (\ref{Fq}) is given by :
\be
F^{h}_{\lambda , q}(x) = F^0_{\lambda , q}(x) h_{ q}(x)
\label{Fx}
\ee
\nn where $h_{ q}(x)$ is a function satisfying
\be
h_{ q}(x) = h_{ q}(qx)
\label{h}
\ee
\nn that is a periodic function of $\ln x/\ln q$ of period 1.
We verify directly that up to a constant the momentums $M(n)$ of $F^0_{\lambda , q}(x)$ are well equal to $\lambda^{-n} q^{-n(n+1)/2}$ for $n\in Z$, as wanted, and this is yet true for $F^{h}_{\lambda , q}(x)$ given in (\ref{Fx}) as implied by the definition of $h_{ q}(x)$.

\nn As the absolute value of a solution of equation (\ref{h}) is also a solution, we can assume that $h_{ q}(x)$ is strictly positive, so that $F^{h}_{\lambda , q}(x)$ is also strictly positive.  

\nn Let us now  consider the integral 
\be
\frac {\int dz d\zz F^{h}_{\lambda , q}(z\zz) f(z) f(\zz)}{\int dz d\zz F^{h}_{\lambda , q}(z\zz )}
\label{int}
\ee
\nn where $f(z)$ is defined in (\ref{f}). Thanks to the positivity of $F^{h}_{\lambda , q}(x)$, we proceed as in Bargmann's paper \cite{bargmann} to prove that the integral (\ref{int}) and the series $\sum_{n\in Z}\mid f_n \mid^2$ are simultaneously divergent or simultaneously convergent to the same value. The function $f(z)$ verifies the inequality :
\be
\mid f(z) \mid \leq C \exp-\left(\frac{\ln^2(x\lambda ^{-1})}{4\ln q}+\frac{\ln(x\lambda ^{-1})}{4}\right)
\sum _{-\infty}^{+\infty} \exp  \frac{\ln q}{2} \left( n+  \frac{1}{2} +\frac{\ln (x\lambda^{-1})}{\ln q }\right)^2 
\label{normf}
\ee
\nn which is essential as implying pointwise convergence from norm convergence.

\nn In particular we use (\ref{normf}) to prove the closedness of the operators
 $z$ and $z^{-1}q^{z\frac{d}{dz}}$. Moreover as we have :
\be
\mid z f(z)\mid ^2 = q \mid z^{-1} q^{z\frac{d}{dz}} f(z) \mid ^2
\ee
\nn we prove that these operators have the same domain of definition and their mutual adjointness  follows easily as in \cite{bargmann}.

\nn So we can summarize the result of this section as follows : for $\psi_{\lambda , q} (N) = \lambda q^{-N}$ the necessary conditions (\ref{moments}) and (\ref{adjoint}) are sufficient conditions for having a Bargmann representation. There exist infinitely many equivalent norms defined by (\ref{F0}) and (\ref{Fx}) for any strictly positive function $h(x)$ verifying (\ref{h})

\section{Generalization of the previous example} 
 
\nn In this section, we point out some directions to extend the results of the previous section
when $\psi (x)$ is of the form :
\be
\psi (x) = \exp \left( \sum_0 ^{2p+1} a_n x^n \right), \quad a_{2p+1}>0
\ee
\nn Then the equation (\ref{mel1}) can be solved and gives :
\be
\hat{F}(\rho ) = \exp \left( \sum_0 ^{2p+1} \frac{a_n}{n+1}B_{n+1}(\rho)
  \right) 
\label{B}
\ee
\nn where $B_{n+1}(\rho)$ are the Bernouilli polynomials.
The term of highest degree in (\ref{B}) is $a_{2p+1} \rho^{2p+2}/(2p+1)$. When $\rho$ is  pure imaginary $\rho = i \sigma , \sigma \in R$, this term is $a_{2p+1}(-1)^{p+1}\sigma^{2p+2}/(2p+1)$, therefore the inverse Mellin transform of $\hat{F}(\rho)$ exists only if $p$ is an even number. The function $F(x)$ thus obtained is always real, but not necessarily positive. Nevertheless, in specific cases, for example when the exponent in (\ref{B}) contains only the term of highest degree, $F(x)$ is actually strictly positive and the Bargmann procedure works as before.

\nn The proof of the closedness of the operators $z$ and $z^{-1}\psi ( z \frac{d}{dz})$ results from (\ref{<}) as before. But as we have 
\be
\mid z f(z) \mid ^2 = \sum_{-\infty}^{\infty} \mid f_n \mid ^2  \psi (n+1)
\ee
\nn and
\be
\mid z^{-1}\psi ( z \frac {d}{dz}) f(z) \mid ^2 = \sum_{-\infty}^{\infty} \mid f_n \mid ^2  \psi (n)
\ee
\nn Since $ \psi(n+1)/\psi(n)$ grows indefinitely as $n \rightarrow \pm \infty$, the domain  $z^{-1}\psi ( z \frac {d}{dz})$ is included in the domain of $z$ but cannot be identical. Nethertheless, the mutual adjointness can be proved as ensured by equation (\ref{adjoint}).

\nn This example points out the following significant fact : apart from its behaviour at infinity, the function $\psi$ must necessarily verify supplementary conditions for admitting a Bargmann representation. Unfortunately, we did not succeed in finding general results in this direction. At least this means that there exist functions $\psi$ such that the equation (\ref{FF}) has no solution fast decreasing at infinity and at the origin simultaneously. 

\section{Deformed algebra associated to a given weight function  }

\nn Conversely, we can use equation (\ref{mel}) to deduce the function $\psi (x)$ from a given $F(x)$. As pointed out just before, we have to prove that the resulting $\psi(x)$ verifies our fundamental assumptions.

\nn We illustrate this construction by the following example.

\nn Let $F(x)$ be 
\be
F(x)= \exp \left(- \nu (\ln x)^{2n} \right), \quad \nu >0, \quad n\geq 1
\ee
\nn The Mellin transform is a strictly positive function given by :
\be
\hat{F}(\rho ) = \int_{-\infty}^{+\infty} \exp \left( -\nu t^{2n} \right) \exp (\rho t) dt
\ee
\nn For $\rho >0$, let us write :
\be
\hat{F}(\rho ) = \int_{-\infty}^{+\infty} \exp \left(- \nu (t+\mu)^{2n} + \rho (t+\mu )\right) dt
\label{hatF}
\ee
\nn with $\mu = \left(\frac{\rho}{2n\nu}\right)^{\frac{1}{2n-1}}$ and take as new variable $u=t \mu ^{n-1}$.
Thus, we get :
\be
\hat{F}(\rho ) = e^{ (2n-1)\nu \left(\frac{\rho}{2n\nu}\right)^{\frac{2n}{2n-1}}}\mu^{1-n}\int_{-\infty}^{+\infty}  e^{  -\nu n(2n-1)u^2 -\nu \sum_{p \geq 3} C^{p}_{2n}\frac{u^p}{\mu^{n(p-2)}}} du
\ee
\nn where $C^p_{2n}$ are the binomial coefficients. It results the following asymptotic expansion for $\hat{F}(\rho )$:
\be
\hat{F}(\rho )= \exp\left( (2n-1)\nu \left(\frac{\rho}{2n\nu}\right)^{\frac{2n}{2n-1}}\right) 
\rho ^{-\frac{n-1}{2n-1}}\sum_{p\geq 0} a_p \rho^{-p}
\ee
Therefore thanks to (\ref{mel1}),for $x \rightarrow +\infty$, $\psi (x)$ has the asymptotic expansion :
\be
\psi(x )= \exp \left( \frac{x}{2n\nu}\right)^{\frac{1}{2n-1}} \sum_{p\geq 0} b_p x^{-p}
\ee
\nn implying $\lim_{x\rightarrow +\infty} \psi(x)=\infty$.
\nn Now from $\hat{F}(\rho)=\hat{F}(-\rho)$, we deduce $\psi (-x) = \psi(x-1)^{-1}$ and thus  $\lim_ {x\rightarrow -\infty} \psi(x)=0$.

\nn Now we use the functional equation (\ref{eqG}) for getting a qualitative evaluation of the asymptotic behaviour of $G(x)$ defined in (\ref{G}) and directly related to the reproducing kernel of the Bargmann representation. Indeed, starting with $\psi(x\frac{d}{dx}) =\hat{F}(x\frac{d}{dx} +1)/\hat{F}(x\frac{d}{dx})$ we get, after replacing $\hat{F}(\rho )$ by the expression (\ref{hatF}) and after performing the needed dilations :
\be
x\int_{-\infty}^{+\infty} \exp (-\nu t^{2n}+t) G(\exp(t)x) dt = \int_{-\infty}^{+\infty} \exp (-\nu t^{2n}+t) G(\exp(t)x) dt 
\ee
\nn This implies that the integral involved in the previous formula is necessarily divergent and  that $G(x)$ grows faster than $\exp (\nu \ln(x)^{2n}-\alpha \ln (x))$ with $\alpha \geq -2$.
\nn Thus, in this specific example, for a given $F(x)$, we have found the corresponding   deformed algebra for which exists a Bargmann representation.

\section{ Bargmann representations corresponding to different $\psi$ }
\nn In this section we exhibit a relation between the functions $F$ associated to different functions $\psi$. The relatively simple results can be used to derive various functions $F$ corresponding to  various function $\psi$ from a given $F$.

\nn Let $\psi _i , i=1,2$, be two positive functions for which exist Bargmann representations and let ${\cal H}_i, i=1,2$, be the corresponding Hilbert spaces. To any sequence $f_n, n \in Z$ such that $\sum_n \mid f_n \mid ^2 < \infty $ is associated two vectors $f_i(z) \in {\cal H}_i,i=1,2$ :
\be
f_i(z) = \sum_{n \geq 0} f_n z^n \left( \psi_i (n)!\right)^{-1/2} + \sum_{n < 0} f_n z^n \left( \psi_i (n)!\right)^{1/2}
\ee
\nn Thus, we can define a unitary mapping $A_{12}$ from ${\cal H}_2$ on ${\cal H}_1$ by :
\be
A_{12}f_2(z) = f_1(z)
\label{map}
\ee
\nn Specifying this equation to the basis elements we deduce easily that $A_{12}$ can be written :
\be
A_{12} = a_{12} (z\frac{d}{dz})
\ee
\nn where the function $a_{12}(\rho )$ verifies :
\be
a_{12}(\rho) = \left(\frac{\psi_2(\rho )}{\psi_1(\rho)}\right)^{1/2} a_{12}(\rho -1)
\label{A}
\ee
\nn Now, we have, for  $f_1, g_1 \in {\cal H}_1$ :
\be
\begin{array}{lll}
(f_1,g_1)_1 &= \int F_1( z\zz)f_1(z) \overline{g_1(z)} dz d\zz&\\
&= \int F_1( z\zz)a_{12}(z \frac{d}{dz})f_2(z) \overline{a_{12}(z \frac{d}{dz})g_2(z)} dz d\zz&=(f_2,g_2)_2
\end{array}
\ee
\nn according to (\ref{map}). Taking into account the behaviour of $F_1(z\zz)$ at infinity and near the origin, the last integral may be written :
\be
\int \left\{ a_{12}(-\frac{d}{dz}z)a_{12}(-\frac{d}{d\zz}\zz)F_1(z\zz)\right\} f_2(z) \overline{g_2(z)} dz d\zz
\ee
\nn so we get 
 \be
F_2(x)= a_{12}^2 (-1-x\frac{d}{dx}) F_1(x)
\label{F2}
\ee
Conversely, this formula can be used to determine a new function $F_2(x)$ from a known one $F_1(x)$ by giving a priori the function $a_{12}^2(\rho)$. In this case, we get $F_2$ by equation (\ref{F2}) and $\psi _2$ by equation (\ref{A}) :
\be
\psi_2 (\rho) = \psi_1(\rho) \frac{a_{12}^2(\rho)}{a_{12}^2(\rho-1)}
\ee
\nn Obviously, we must be careful in choosing $a_{12}^2(\rho)$ to be sure that the resulting $\psi_2 (\rho)$ fits the conditions for admitting a Bargmann representation. For instance, we choose
\be
a_{12}^2(\rho) = \sum _0^p a_n \exp(\alpha _n \rho ), \quad a_n > 0
\ee
\nn where the $\alpha _n $ are in increasing order.
Then we have :
\be
F_2 (x) = \sum _0^p a_n \exp(\alpha _n) F_1 (x\exp(-\alpha _n))
\ee
\nn and
\be
\psi_2(\rho ) = \psi_1(\rho ) \frac{\sum _0^p a_n \exp(\alpha _n \rho )}{\sum _0^p a_n \exp(\alpha _n (\rho -1))}
\ee
\nn on which appears clearly that the conditions for a Bargmann representation are satisfied.

\section{The ring case}

\nn Let us consider now the case where the domain of convergence of the coherent states is a ring $\cal{R}$ in the complex plane defined by $r_1 < \mid z \mid < r_2$. We are looking for $F(z\zz )$ such that 
\be
\int_{\cal{R}} F(z\zz ) \mid \zz>< \zz \mid dz d\zz = 1
\ee
\nn We can always suppose that $F(z\zz )$ is identically zero outside $\cal{R}$. Then, we can develop the same consideration as above and from the mutual adjointness of $z$ and $z^{-1}\psi (z\frac{d}{dz})$, we get the same fuctional equation (\ref{FF}) with the limit conditions that $F(x)$ is zero outside the interval $r_{1}^2 < x <r_{2}^2$.

\nn Let us illustrate the situation with two examples :

a) $\psi (x) = 1 + q^x , \quad q>1$

\nn In this case $\cal{R}$ is the  exterior of the disk of radius $1$ and the functional equation (\ref{FF}) reads :
\be
(qx-1) F(qx) = F(x)
\label{qf}
\ee
\nn So, if $F(x) \equiv 0$ for $x< 1$, then $F(x) \equiv 0$ for $x<q$. 
Using (\ref{qf}), we obtain 
\be
(q^{n+1}x-1) F(q^{n+1}x) = F(q^n x)
\ee
\nn With this equation, we prove recursively that $F(x)$ is identically zero in the complex plane.

b) $\psi (x) = \frac{q^x}{1+ q^x}, \quad q>1$

\nn In this case $\cal{R}$ is the disk of radius $1$. The functional equation (\ref{FF}) reads :
\be
x F(x) = (1- q^{-1} x) F(q^{-1}x)
\ee
\nn If $F(x) \equiv 0$ for $ x>1$, then $F(x) \equiv 0$ for $x> q^{-1}$ and as previously we show recursively that $F(x)$ is identically zero in the complex plane.

\nn Although we have not obtained a general proof, we conjecture that the results proved in the two previous examples are general : when the coherent states are restricted to a ring of the complex plane, we cannot obtain a Bargmann representation because the function involved in the resolution of the identity is identically equal to zero. 
We think that the phenomenon is general thanks to the fact that the non-local operator $\psi(-x\frac{d}{dx})$  does not leave invariant the limit conditions. One can also think that this negative result is depending on our restricted definition for Bargmann representation, and that we would have advantage to generalize it, by introducing another type of integration like for instance the q-integration used in \cite{gray},\cite{bracken} that precisely is the inverse of a non-local operator.

\section{Conclusion}
\nn The deformed oscillators considered in this paper are essentially different from the well-known harmonic oscillator, since their irreducible representations never admit a fundamental vector. Nevertheless it has been possible for some of these exotic cases to realize their representations in a Hilbert space of holomorphic functions in complete analogy with the well-known Bargmann construction, taken in the more restricted sense i.e. by excluding any new concept of integration as in \cite{gray},\cite{bracken}. We have given necessary conditions to be verified by the weight function, or equivalently, by its Mellin transform. These conditions were shown to be sufficient in several specific cases, including the so called q-oscillator. Meanwhile, we encountered cases where we proved the non-existence of any weight function. Although our results have some degree of generality, nevertheless we did not get a general statement giving a complete characterization of the deformed alg!
 eb!
ras, considered in this paper, admitting a Bargmann representation.

\end{document}